\documentclass[12pt,leqno,draft]{article}

\newtheorem{theorem}{Theorem}
\newtheorem{lemma}[theorem]{Lemma}
\newtheorem{proposition}[theorem]{Proposition}
\newtheorem{definition}[theorem]{Definition}
\newtheorem{corollary}[theorem]{Corollary}

\newcommand{\begintheorem}{\addtocounter{equation}{1}\begin{theorem}}
\newcommand{\beginlemma}{\addtocounter{equation}{1}\begin{lemma}}
\newcommand{\beginproposition}{\addtocounter{equation}{1}\begin{proposition}}
\newcommand{\begindefinition}{\addtocounter{equation}{1}\begin{definition}}
\newcommand{\begincorollary}{\addtocounter{equation}{1}\begin{corollary}}

\begin{document}

\title{Some notes about matrices, 2}

\author{Stephen William Semmes 	\\
	Rice University		\\
	Houston, Texas}

\date{}

\maketitle

\tableofcontents

\section{Lattices in ${\bf R}^n$}
\label{section on lattices in R^n}
\setcounter{equation}{0}

	As usual, ${\bf R}$ denotes the real numbers, ${\bf Z}$
denotes the integers, and ${\bf R}^n$, ${\bf Z}^n$ consist of
$n$-tuples of real numbers and integers, respectively.  Sometimes we
might refer to ${\bf Z}^n$ as the \emph{standard integer lattice} in
${\bf R}^n$.  If we say that $L$ is a \emph{lattice} in ${\bf R}^n$,
then we mean that there is an invertible linear transformation $A$ on
${\bf R}^n$ such that
\begin{equation}
	L = A({\bf Z}^n).
\end{equation}

	If $L$ is a lattice in ${\bf R}^n$, then we can form the
quotient ${\bf R}^n / L$.  That is, two vectors $x$, $y$ in ${\bf
R}^n$ are identified in the quotient if their difference $x - y$ lies
in $L$.  In particular, we get a canonical quotient mapping
\begin{equation}
	\rho : {\bf R}^n \to {\bf R}^n / L
\end{equation}
which sends a vector $x$ in ${\bf R}^n$ to the corresponding
element of the quotient.

	Now, with respect to ordinary vector addition, ${\bf R}^n$ is
an abelian group, and a lattice $L$ is a subgroup of ${\bf R}^n$.  We
can think of the quotient ${\bf R}^n / L$ as a quotient in the sense
of group theory.  The quotient is an abelian group under addition, and
the canonical quotient mapping is a group homomorphism.

	We can also look at the quotient ${\bf R}^n / L$ in terms of
topology.  Namely, it inherits a topology from the one on ${\bf R}^n$
so that the canonical quotient mapping is an open continuous mapping,
which means that both images and inverse images of open sets are open
sets, and indeed the canonical quotient mapping is a nice covering
mapping, so that for every point $x$ in ${\bf R}^n$ there is a
neighborhood $U$ of $x$ in ${\bf R}^n$ such that the restriction of
$\rho$ to $U$ is a homeomorphism from $U$ onto the open set $\rho(U)$ in
${\bf R}^n / L$.  For that matter we can think of ${\bf R}^n / L$ as a
smooth manifold, with the quotient mapping $\rho$ as a smooth mapping
which is a local diffeomorphism.

	Suppose that $L_1$, $L_2$ are lattices in ${\bf R}^n$, and
let
\begin{equation}
	\rho_1 : {\bf R}^n \to {\bf R}^n / L_1, \quad
		\rho_2 : {\bf R}^n \to {\bf R}^n / L_2
\end{equation}
be the corresponding canonical quotient mappings.  If $A$ is an
invertible linear transformation on ${\bf R}^n$ such that
\begin{equation}
	A(L_1) = L_2,
\end{equation}
then we get an induced mapping
\begin{equation}
	\widehat{A} : {\bf R}^n / L_1 \to {\bf R}^n / L_2.
\end{equation}
This mapping is a group isomorphism and a homeomorphism, and even
a diffeomorphism, which satisfies the obvious compatibility condition
with the corresponding canonical quotient mappings $\rho_1$, $\rho_2$,
namely $\rho_1 \circ A = \widehat{A} \circ \rho_2$.

	When $n = 1$, one can consider the lattice $2 \pi {\bf Z}$
consisting of integer multiples of $2 \pi$, and it is customary to
identify ${\bf R} / 2 \pi {\bf Z}$ with the unit circle ${\bf T}$ in
the complex numbers ${\bf C}$,
\begin{equation}
	{\bf T} = \{z \in {\bf C} : |z| = 1\},
\end{equation}
where $|z|$ denotes the usual modulus of $z \in {\bf C}$, $|z| = (x^2
+ y^2)^{1/2}$ when $z = x + i \, y$, $x, y \in {\bf R}$.  More
precisely, $\exp (i \, t)$ is an explicit version of the canonical
quotient mapping from ${\bf R} / 2 \pi {\bf Z}$ onto ${\bf T}$ with
respect to this identification, which is a local diffeomorphism and a
group homomorphism using the group structure of multiplication on
${\bf T}$.  In general, we can identify ${\bf R}^n / 2 \pi {\bf Z}^n$
with ${\bf T}^n$, the $n$-fold Cartesian product of ${\bf T}$, where
$2 \pi {\bf Z}^n$ denotes the lattice of points whose coordinates are
all integer multiples of $2 \pi$.

	Suppose that $L$ is a lattice in ${\bf R}^n$.  Also let $A$ be
an invertible linear mapping on ${\bf R}^n$ such that $A(2 \pi {\bf
Z}^n) = L$.  Thus $\widehat{A}$ is a group isomorphism and a
diffeomorphism from ${\bf R}^n / 2 \pi {\bf Z}^n \cong {\bf T}^n$ onto
${\bf R}^n / L$.

	There is a more precise way to look at the quotient of ${\bf
R}^n$ by a lattice, which is to say that the quotient space has a kind
of local affine structure.  That is, there is a local affine structure
in which the canonical quotient mapping is considered to be locally
affine, and which permits one to say when a curve in the quotient is
locally a straight line segment, like an arc on a line, and when it
has locally constant speed, etc.  If $L_1$, $L_2$ are lattices in
${\bf R}^n$ and $A$ is an invertible linear mapping on ${\bf R}^n$
such that $A(L_1) = L_2$, then the induced mapping $\widehat{A}$ from
${\bf R}^n / L_1$ onto ${\bf R}^n / L_2$ preserves this local affine
structure on the quotient spaces.

	There is an even more precise way to look at the quotient
${\bf R}^n / L$ of ${\bf R}^n$ by a lattice $L$, which is that it has
a local flat geometric structure, induced from the one on ${\bf R}^n$.
With respect to this structure one can make local measurements of
lengths, volumes, and angles, like the length of a curve, the angle at
which two curves meet at a point, or the volume of a nice subset.  In
technical terms this can be seen as a special case of a
\emph{Riemannian metric}.

	In particular, one can define the volume of such a quotient
${\bf R}^n / L$, where the volume of ${\bf R}^n / {\bf Z}^n$ is equal
to $1$, and the volume of ${\bf R}^n / 2 \pi {\bf Z}^n$ is equal to
$(2 \pi)^n$.  In general, if $L_1$, $L_2$ are lattices in ${\bf R}^n$
and $A$ is an invertible linear transformation on ${\bf R}^n$ such
that $A(L_1) = L_2$, then the volume of ${\bf R}^n / L_2$ is equal to
$|\det A|$ times the volume of ${\bf R}^n / L_1$, and more generally
if $E$ is a nice subset of ${\bf R}^n / L_1$, then the volume of
$\widehat{A}(E)$ in ${\bf R}^n / L_2$ is equal to $|\det A|$ times the
volume of $A$ in ${\bf R}^n / L_1$.  This is a variant of the fact
that on ${\bf R}^n$ a linear transformation $A$ distorts volumes by
a factor of $|\det A|$, where $\det A$ denotes the determinant of $A$.

	Suppose that $L_1$, $L_2$ are lattices in ${\bf R}^n$, and
that $T$ is a linear transformation on ${\bf R}^n$ such that $T(L_1) =
L_2$.  Recall that $T$ is an \emph{orthogonal transformation} on ${\bf
R}^n$ if $T$ is invertible with inverse given by the adjoint, also
known as the transpose, of $T$, and that this is equivalent to saying
that $T$ preserves the standard norm of vectors in ${\bf R}^n$, and
the standard inner product of vectors in ${\bf R}^n$.  In other words,
orthogonal transformations on ${\bf R}^n$ are linear mappings which
preserve the geometry in ${\bf R}^n$, and for the lattices $L_1$, $L_2$
and the quotients of ${\bf R}^n$ by them we have that the induced
mapping $\widehat{T}$ from ${\bf R}^n / L_1$ onto ${\bf R}^n / L_2$
preserves the geometry as well.

	In short, quotients of ${\bf R}^n$ by lattices are the same in
terms of group structure, topological and even smooth structure, and
affine structure, and not in general for more precise geometry.  The
volume of the quotient space is one basic parameter that one can
consider.  It is also interesting to look at closed curves in the
quotient which are locally flat, their lengths, the angles at which
they meet, and so on.

\section{Lattices in ${\bf C}^n$}
\label{section on lattices in C^n}
\setcounter{equation}{0}

	As before we write ${\bf C}^n$ for the $n$-tuples of real
numbers.  We can identify ${\bf C}^n$ with ${\bf R}^{2n}$ in the usual
manner, so that the real and imaginary parts of the $n$ components of
an element of ${\bf C}^n$ give rise to the $2n$ components of an
element of ${\bf R}^{2n}$.  By a lattice in ${\bf C}^n$ we mean a
lattice in ${\bf R}^{2n}$ which is then identified with ${\bf C}^n$.

	Let us write ${\bf Z}[i]$ for the \emph{Gaussian integers},
which are complex numbers of the form $a + i \, b$, where $a$, $b$ are
integers.  We also write $({\bf Z}[i])^n$ for the lattice in ${\bf
C}^n$ consisting of $n$-tuples of Gaussian integers.  We call this the
standard integer lattice in ${\bf C}^n$.

	If $L$ is a lattice in ${\bf C}^n$, then the quotient ${\bf
C}^n / L$ inherits a complex structure from ${\bf C}^n$.  If $L_1$,
$L_2$ are lattices in ${\bf C}^n$ and $A$ is an invertible
complex-linear transformation on ${\bf C}^n$ such that $A(L_1) = L_2$,
then $A$ induces a mapping $\widehat{A}$ from ${\bf C}^n / L_1$ to
${\bf C}^n / L_2$ which preserves this complex structure.  Thus,
although lattices in ${\bf C}^n$ can be defined as images of the
standard integer lattice in ${\bf C}^n$ under invertible real linear
transformations, complex structures and complex linear transformations
play an important role.

	Recall that a unitary transformation on ${\bf C}^n$ is an
invertible complex-linear transformation whose inverse is given by the
adjoint.  This is the same as saying that the linear transformation
preserves the standard Hermitian inner product on ${\bf C}^n$, and it
is also equivalent to saying that the transformation is both
complex-linear and an orthogonal transformation with respect to the
standard identification of ${\bf C}^n$ with ${\bf R}^{2n}$.  This
leads to another and more precise relationship between lattices in
${\bf C}^n$ and their quotients, i.e., having a unitary transformation
which takes one lattice to another, and which then induces a nice
transformation between the corresponding quotients.

\section{Spaces of lattices}
\label{section on spaces of lattices}
\setcounter{equation}{0}

	Let us write $GL({\bf R}^n)$ for the \emph{general linear
group} on ${\bf R}^n$, consisting of the invertible linear
transformations on ${\bf R}^n$, with composition as the group
structure.  The \emph{special linear group} $SL({\bf R}^n)$ is the
subgroup of $GL({\bf R}^n)$ of linear transformations with determinant
equal to $1$.  The \emph{orthogonal group} $O({\bf R}^n)$ is the
subgroup of $GL({\bf R}^n)$ of orthogonal linear transformations on
${\bf R}^n$, and the \emph{special orthogonal group} $SO({\bf R}^n)$
consists of orthogonal linear transformations whose determinant is
equal to $1$.

	Consider the quotient space $O({\bf R}^n) \backslash GL({\bf
R}^n)$, in which two invertible linear transformations on ${\bf R}^n$
are identified if one can be written as an orthogonal linear
transformation times the other.  We can identify this quotient space
with the space of symmetric linear transformations on ${\bf R}^n$
which are positive definite, through the mapping
\begin{equation}
	T \mapsto T^* \, T.
\end{equation}
In other words, if $T$ is an invertible linear transformation on ${\bf
R}^n$, then $T^* \, T$ is a symmetric linear transformation on ${\bf
R}^n$ which is positive-definite, $T_1^* \, T_1 = T_2^* \, T_2$ for
$T_1, T_2 \in GL({\bf R}^n)$ if and only if $T_2 = R \, T_1$ for some
orthogonal transformation $R$, and every symmetric linear transformation
on ${\bf R}^n$ which is positive-definite can be expressed as $T^* \, T$
for an invertible linear transformation $T$.

	Similarly, the quotient $SO({\bf R}^n) \backslash SL({\bf
R}^n)$ can be identified with the space $\mathcal{M}({\bf R}^n)$ of
symmetric linear transformations on ${\bf R}^n$ which are positive
definite and have determinant equal to $1$.  Let us also write
$\Sigma({\bf R}^n)$ for the elements of $SL({\bf R}^n)$ whose matrices
with respect to the standard basis have integer entries.  The inverse
of a linear transformation in $\Sigma({\bf R}^n)$ also lies in
$\Sigma({\bf R}^n)$, because Cramer's rule gives a formula for the
matrix of the inverse which shows that it has integer entries when
the original matrix has integer entries and determinant equal to $1$.

	Elements of $\Sigma({\bf R}^n)$ can be described as the
invertible linear transformations which take ${\bf Z}^n$ onto itself.
The quotient $SL({\bf R}^n) / \Sigma({\bf R}^n)$ describes the space
of lattices $L$ in ${\bf R}^n$ such that the corresponding quotient
${\bf R}^n / L$ has volume equal to $1$ and for which there is an
extra piece of data concerning orientation, and the double quotient
$SO({\bf R}^n) \backslash SL({\bf R}^n) / \Sigma({\bf R}^n)$ deals
with these lattices up to equivalence under rotation.  Of course there
are a lot of variants of these themes, which may involve complex
structures in particular.

\end{document}